\documentclass[12pt,amsfonts]{article}
\usepackage[margin=1in]{geometry}  
\usepackage{graphicx}              
\usepackage{amsmath}               
\usepackage{amsfonts}              
\usepackage{amsthm}                
\newtheorem{claim}{Claim}[section]
\newtheorem{thm}{Theorem}[section]

\newtheorem{lem}{Lemma}[section]

\newtheorem{cor}{Corollary}[section]
\newcommand{\fim}{\hfill\rule{2mm}{2mm}}

\numberwithin{equation}{section}

\begin{document}

\title{Existence of positive solution for a nonlinear elliptic equation with saddle-like potential and nonlinearity with exponential critical growth in $\mathbb{R}^{2}$. }
\author{\sf Claudianor O. Alves\thanks{Research of C. O. Alves partially supported by  CNPq 304036/2013-7  and INCT-MAT} 
\\
\small{Universidade Federal de Campina Grande, }\\
\small{Unidade Acad\^emica de Matem\'atica } \\
\small{CEP:58429-900, Campina Grande-PB, Brazil}\\}
\date{}
\maketitle
\begin{abstract}
In this paper, we use variational methods to prove the existence of positive solution for the following class of elliptic equation 
$$
-\epsilon^{2}\Delta{u}+V(z)u=f(u) \,\,\, \mbox{in} \,\,\, \mathbb{R}^{2},
$$
where $\epsilon >0$ is a positive parameter, $V$ is a saddle-like potential and $f$ has an exponential critical growth.
\end{abstract}

\vspace{0.5 cm}


\section{Introduction}

Many recent studies have focused on the nonlinear Schr\"{o}dinger
equation
$$
i\epsilon \displaystyle \frac{\partial \Psi}{\partial t}=-\epsilon ^{2}\Delta
\Psi+(V(z)+E)\Psi-f(\Psi)\,\,\, \mbox{for all}\,\,\, z \in
\Omega,\eqno{(NLS)}
$$
where $\epsilon > 0$, $\Omega$ is a domain in $\mathbb{R}^N$, not
necessarily bounded, with empty or smooth boundary. Knowledge of the
solutions for the elliptic equation
$$
\ \  \left\{
\begin{array}{l}
- \epsilon^{2} \Delta{u} + V(z)u=f(u)
\ \ \mbox{in} \ \ \Omega,
\\
u = 0 \ \ \mbox{on} \ \ \partial\Omega,
\end{array}
\right.
\eqno{(P)}
$$
has a great importance in the study of standing-wave solutions
of (NLS). In recent years, the existence and concentration of
positive solutions for general semilinear elliptic equations
$(P_\epsilon)$ for the case $N \geq 2$ have been extensively
studied, see for example, Bartsch, Pankov, and Wang \cite{BPW},
Bartsch and Wang \cite{BW0},   Floer and Weinstein \cite{FW}, Oh
\cite{Oh1}, Rabinowitz \cite{R}, Wang \cite{W}, Alves and Souto
\cite{Alves10}, del Pino and Felmer \cite{Pino}, del Pino, Felmer and Miyagaki \cite{PFM}, Alves, do \'{O} and
Souto \cite{AOS}, do \'O and Souto \cite{OS}, Alves and Soares \cite{ASS1} and their references.

In \cite{R}, by a mountain pass argument,  Rabinowitz proved the
existence of positive solutions of $(NLS)$, for $\epsilon > 0$
small, whenever
\[
\liminf_{|z| \rightarrow \infty} V(z) > \inf_{z \in
	\mathbb{R}^N}V(z)=\gamma >0.
\]
Later Wang \cite{W} showed that these solutions concentrate at
global minimum points of $V$ as  $\epsilon$ tends to 0.

In \cite{Pino}, del Pino and Felmer established  the existence of positive  solutions which concentrate around local minimum of $V$ by introducing a penalization method. More precisely, they assumed that there is an open and bounded set $\Lambda$ compactly contained  in $\Omega$ satisfying
$$
0< \gamma \leq V_0 =\inf_{z\in\Lambda}V(z)< \min_{z \in
	\partial\Lambda}V(z). \eqno{(V_{1})}
$$

In \cite{PFM}, del Pino, Felmer and Miyagaki considered the case where potential $V$ has a geometry like saddle, essentially they assumed the following conditions on $V$: First of all, they fixed two subspaces $X,Y  \subset \mathbb{R}^{N}$ such that
$$
\mathbb{R}^{N}=X \oplus Y.
$$  
By supposing that $V$ is bounded,  they fixed $c_0,c_1>0$  satisfying \\
$$
\displaystyle c_0=\inf_{z \in \mathbb{R}^{N}}V(z)>0
$$
and
$$
c_1=\displaystyle \sup_{x \in X}V(x).
$$
Furthermore, they also supposed that $V \in C^{2}(\mathbb{R}^{N})$ and it  verifies the following geometry conditions:  \\

\noindent $(V_1)$  
$$
c_0=\inf_{R>0}\sup_{x \in \partial B_R(0) \cap X }V(x)<\inf_{y \in Y}V(y).
$$
\noindent $(V_2)$ \quad The functions $V, \frac{\partial V}{\partial x_i}$ and $\frac{\partial^{2} V}{\partial x_i \partial x_j}$ are bounded in $\mathbb{R}^{N}$ for all $i,j \in \{1,...,N\}$.  \\

\noindent $(V_3)$ \quad $V$ satisfies the Palais-Smale condition, that is, if $(x_n) \subset \mathbb{R}^{N}$ is a sequence such that $(V(x_n))$ is bounded and $\nabla V(x_n) \to 0$, then $(x_n)$ possesses a convergent subsequence in $\mathbb{R}^{N}$. \\

Using the above conditions on $V$, and supposing that 
$$
c_1<2^{\frac{2(p-1)}{N+2-p(N-2)}}c_0,
$$
del Pino, Felmer and Miyagaki showed  the existence of positive solutions for the following problem  
$$
- \epsilon^{2} \Delta{u} + V(z)u=|u|^{p-2}u
\ \ \mbox{in} \ \ \mathbb{R}^{N},
$$
where $p \in (2,2^{*})$ if $N \geq 3$ and $p \in (2,+\infty)$ if $N=1,2$, for $\epsilon>0$ small enough. The main tool used was the variational method, more precisely, the authors found critical point of the functional
$$
E_{\epsilon}(u)=\int_{\mathbb{R}^{N}}(\epsilon^{2}|\nabla u|^{2}+|u|^{2})\,dx
$$ 
on the manifold
$$
\mathcal{M} =\left\{u \in H^{1}(\mathbb{R}^{N}) \cap P \,:\,\int_{\mathbb{R}^{N}}|u|^{p}\,dx=1 \right\},
$$
where $P$ denotes the cone of nonnegative functions of $H^{1}(\mathbb{R}^{N})$. 

Motivated by the results obtained in \cite{PFM},  our goal is to show the existence of positive solution for the following class of elliptic equation
$$
-\epsilon^{2}\Delta{u}+V(x)u=f(u) \,\,\, \mbox{in} \,\,\, \mathbb{R}^{2}, \eqno{(P_\epsilon)}
$$
where $\epsilon >0$ is a positive parameter and $V,f$ are functions verifying some conditions.  Here, the potential $V$ has the same geometry considered in \cite{PFM}, that is, $V$ satisfies the  conditions $(V_1)-(V_3)$. Related to function $f$, we will suppose that it has an exponential critical growth. 

In  $\mathbb{R}^2$, the natural growth restriction on the function $f$  is given by the inequality
of Trudinger and Moser \cite{M,T}. More precisely, we say that a
function $f$ has an exponential critical growth if there is $\alpha_0 >0$ such that
$$ \lim_{|s| \to \infty} \frac{|f(s)|}{e^{\alpha s^{2}}}=0
\,\,\, \forall\, \alpha > \alpha_{0}\quad \mbox{and} \quad
\lim_{|s| \to \infty} \frac{|f(s)|}{e^{\alpha s^{2}}}=+ \infty
\,\,\, \forall\, \alpha < \alpha_{0}.
$$
We would like to mention that problems involving exponential critical growth have received a special attention at last years, see for example, \cite{A, AdoOM,Cao,doORuf,DMR,DdOR} for semilinear elliptic equations, and \cite{1,2,3,4,5} for quasilinear equations. 

Hereafter, $f \in C^{1}(\mathbb{R})$ and it satisfies the following conditions:

\begin{enumerate}
\item[$(f_1)$] There is $C>0$ such that
$$
|f(s)|,|f'(s)| \leq Ce^{4\pi |s|^2}\ \ \mbox{for all}\ \ s\in \mathbb{R}.
$$
\item[$(f_2)$] $\displaystyle \lim_{s\rightarrow
0}\dfrac{f(s)}{s}=0$.
\item[$(f_3)$] There is $\theta>2$ such that
$$
0<\theta F(s):=\theta\int_0^{s}f(t)dt\leq sf(s),\ \ \mbox{for all}\ \ s \in \mathbb{R} \setminus \{0\}. 
$$
\item[$(f_4)$] The function $s\rightarrow\dfrac{f(s)}{s}$ is
strictly increasing in $|s|>0$. 
\item[$(f_5)$] There exist constants $p>2$ and $C_p>0$ such  that
$$
f(s)\geq C_p s^{p-1}\ \ \mbox{for all}\ \ s >0,
$$
where
$$
C_p>\left[8 \beta_p\left(\frac{2\theta}{\theta-2}\right)\frac{1}{\min\{1,c_0\}}\right]^{(p-2)/2},
$$
with
$$
\beta_p=\inf_{\mathcal{N}_\infty}J_\infty,
$$
$$
\mathcal{N}_\infty=\{u\in H^1(\mathbb{R}^{2}) \setminus \{0\}:\ J_{\infty}'(u)u=0\}
$$
and
$$
J_{\infty}(u)=\dfrac{1}{2}\int_{\mathbb{R}^{2}}\left(|\nabla u|^2+|V|_\infty|u|^2\right)dx-\dfrac{1}{p}\int_{\mathbb{R}^{2}}|u|^pdx .
$$

\end{enumerate}

 \vspace{0.5 cm}

Before to state our main result, we need to fix some notations. Hereafter, if $A \in \mathbb{R}$, we denote by $J_{A}$ the functional given by
$$
J_{A}(u)=\frac{1}{2}\int_{\mathbb{R}^{2}}\left(|\nabla u|^{2}+A|u|^{2}\right)\,dx-\int_{\mathbb{R}^{2}}F(u)\,dx,
$$
which are defined in $H^{1}(\mathbb{R}^{2})$. Moreover, let us denote by $m(A)$ the mountain pas level associated with $J_A$, which possesses the following characterizations   
$$
m(A)=\inf_{u \in H^{1}(\mathbb{R}^{2}) \setminus \{0\}}\left\{\max_{t \geq 0}J_A(tu) \right\}=\inf_{u \in \mathcal{M}_A}J_A(u),
$$
where $\mathcal{M}_A$ is the Nehari Manifold associated with $J_A$, given by
$$
\mathcal{M}_A=\{u \in H^{1}(\mathbb{R}^{2})\setminus\{0\}\,:\,J'_A(u)u=0\}.
$$

The main result in the present paper is the following

\begin{thm} \label{T1} Assume that $(V_1)-(V_3)$ and $(f_1)-(f_5)$ hold. If 
$$ 
m(V(0)) \geq 2m(c_0) \quad \mbox{and} \quad c_1 \leq \left[1 + \frac{3}{5}\left(\frac{1}{2}-\frac{1}{\theta}\right)\right]c_0, \leqno{(V_4)} 
$$
then there is $\epsilon_0>0$ such that $(P_\epsilon)$ has a positive solution for all $\epsilon \in (0,\epsilon_0]$.  
\end{thm}

The  inspiration to prove the Theorem \ref{T1} comes from \cite{PFM}, however  it is important to say that we are working with exponential critical growth, then our estimates for this class of problem are very delicate and different from those used in the above mentioned paper. Here, we have proved a lot of estimates that do not appear in \cite{PFM}, for more details see Section 2. In Section 3, we minimize the energy function on the  Nehari manifold, and to this end, we modify the barycenter function of a way more convenient for our problem, see Section 3.

\vspace{0.5 cm}

Before to conclude this introduction, we would like point out that using the change variable $v(x)=u(\epsilon x)$, it is possible to prove that $(P_\epsilon)$ is equivalent to the following problem
$$
-\Delta{u}+V(\epsilon x)u=f(u) \,\,\, \mbox{in} \,\,\, \mathbb{R}^{2}. \eqno{(P_\epsilon)'}
$$
In what follows, we denote by $I_\epsilon$ the energy functional associated with $(P_\epsilon)'$, that is, 
$$
I_\epsilon(u)=\frac{1}{2}\int_{\mathbb{R}^{2}}(|\nabla u|^{2}+V(\epsilon x)|u|^{2})\,dx-\int_{\mathbb{R}^{2}}F(u)\,dx, \quad \forall u \in H^{1}(\mathbb{R}^{2}).
$$
This way, $u \in H^{1}(\mathbb{R}^{2})$ is a weak solution for $(P_\epsilon)'$ if, and only if, $u$ is a critical point for $I_\epsilon$.  

\vspace{0.5 cm}

\noindent \textbf{Notation:} In this paper, we use the following
notations:
\begin{itemize}
	\item  The usual norms in $L^{t}(\mathbb{R}^{2})$ and $H^{1}(\mathbb{R}^{2})$ will be denoted by
	$|\,. \,|_{t}$ and $\|\;\;\;\|$ respectively.

	\item   $C$ denotes (possible different) any positive constant.
	
	\item   $B_{R}(z)$ denotes the open ball with center at $z$ and
	radius $R$.

\end{itemize}

\section{Technical results}

 In this section, we will prove some technical lemmas, which are crucial in our approach. Since we will work with exponential critical growth in whole $\mathbb{R}^{2}$,
 a key inequality in our arguments is the Trudinger-Moser inequality in whole space $\mathbb{R}^2$ due to Cao \cite{Cao}, which has the following statement
 
 \begin{lem} \label{Cao} For all $u \in H^1(\mathbb{R}^2)$ and $\alpha>0$, 
 \begin{equation} \label{X2}
 \int_{\mathbb{R}^2}(e^{\alpha|u|^2}-1)\,dx<+\infty. 
 \end{equation}
 Furthermore, if $\alpha \leq 4\pi$ and $|u|_{2}\leq M$, there exists
 a constant \linebreak $C_1=C_1(M,\alpha)$ such that
 \begin{equation} \label{X4}
 \sup_{ |\nabla u|_{2}\leq 1}
 \int_{\mathbb{R}^2}(e^{\alpha|u|^2}-1)\,dx\leq C_1.
 \end{equation}
 \end{lem}

Using the above lemma, we are able to prove some technical lemmas. The first of them is crucial in the study of the $(PS)$ condition for $I_\epsilon$. 
\begin{lem} \label{alphat11} Let $(u_{n})$ be a sequence in $H^{1}(\mathbb{R}^{2})$ with 
	$$
	\limsup_{n \to +\infty} \|u_n \|^{2}  < 1.
	$$ 
	Then, there exist  $t > 1$ with $t \approx 1$ and $C > 0$  satisfying
	\[
	\int_{\mathbb{R}^{N}}\left(e^{\alpha|u_n|^{2}} - 1 \right)^t dx \leq C, \,\,\,\,\forall\, n \in \mathbb{N}. 
	\]

\end{lem}
\noindent {\bf Proof.} \, 
As
\[
\limsup_{n\rightarrow\infty} ||u_n||^2 < 1,
\]
there are $m >0$ and $n_0\in \mathbb{N}$ verifying 
\[
||u_n||^{2} < m < 1,
\,\,\,\,\forall\,  n \geq n_0.
\]
Fix  $t > 1$ with $t \approx 1$ and $ \beta> t$ satisfying $\beta m < 1. $ Then, there exists $C=C(\beta)>0$ such that
$$
\int_{\mathbb{R}^{2}}\left(e^{4 \pi|u_n|^{2}} -1 \right)^t dx   
\leq C\int_{\mathbb{R}^{2}}\left(e^{4 \beta m  \pi(\frac{|u_n|}{||u_n||})^{2}} - 1 \right) dx,
$$
for every  $n \geq n_0$.  Hence, by Lemma \ref{Cao},
$$
\int_{\mathbb{R}^{2}}\left(e^{4 \pi |u_n|^{2}} -1 \right)^t dx \leq C_1 \, \,\,\,\,\forall\,  n \geq n_0,
$$
for some positive constant $C_1$. Now, the lemma follows fixing
$$
C=\max\left\{C_1, \int_{\mathbb{R}^{2}}\left(e^{4 \pi |u_{1}|^{2}} -1 \right)^t dx,....,\int_{\mathbb{R}^{2}}\left(e^{4 \pi |u_{n_0}|^{2}} -1 \right)^t dx \right\}.
$$
\fim

\begin{cor} \label{Convergencia em limitados} Let $(u_{n})$ be a sequence in $H^{1}(\mathbb{R}^{2})$ with 
$$
	\limsup_{n \to +\infty} \|u_n \|^{2}  < 1.
$$ 
If $u_n \rightharpoonup u$ in $H^{1}(\mathbb{R}^{2})$ and $u_n(x) \to u(x)$ a.e in $\mathbb{R}^{2}$, then
$$
F(u_n) \to F(u) \quad \mbox{in} \quad L^{1}(B_R(0)), \quad \forall R>0.
$$ 	
\end{cor}
\noindent {\bf Proof.} By $(f_1)$, for each $\beta >1$, there is $C>0$ such that
$$
|F(s)| \leq C(e^{4 \beta \pi |s|^{2}}-1) \quad \forall s \in \mathbb{R},
$$
from where it follows that, 
\begin{equation} \label{Domina}
|F(u_n)| \leq C(e^{4 \beta \pi |u_n|^{2}}-1), \quad \forall n \in \mathbb{N}.
\end{equation}
Setting
$$
h_n(x)=C(e^{4 \beta \pi |u_n(x)|^{2}}-1),
$$
we can fix $\beta, q>1$ with $\beta,q \approx 1$ such that 
$$
h_n \in L^{q}(\mathbb{R}^{2}) \quad \mbox{and} \quad \sup_{n \in \mathbb{N}}|h_n|_q<+\infty, 
$$
which is an immediate consequence of Lemma \ref{alphat11}. Therefore, for some subequence of $(u_n)$, still denoted by itself, we derive that
$$
h_n \rightharpoonup h=C(e^{4 \beta \pi |u|^{2}}-1) \quad \mbox{in} \quad L^{q}(\mathbb{R}^{2}).
$$
As $h_n, h \geq 0$, the last limit yields 
$$
h_n \to h \quad \mbox{in} \quad L^{1}(B_R(0)), \quad \forall R>0.
$$
Gathering the last limit  with (\ref{Domina}), we get 
$$
F(u_n) \to F(u) \quad \mbox{in} \quad L^{1}(B_R(0)), \quad \forall R>0.
$$

\fim

The next lemma is a Lions type result for exponential critical growth due to Alves, do \'O and Miyagaki  \cite{AdoOM}.

\begin{lem} \label{lions} Let $(u_n) \subset H^{1}(\mathbb{R}^{2})$ be a sequence  with 
$$
\limsup_{n \to +\infty} \|u_n \|^{2}  < 1.
$$
If there is $R>0$ such that
$$
\lim_{n \to +\infty}\sup_{z \in \mathbb{R}^{2}}\int_{B_R(z)}|u_n|^{2}\,dx=0,
$$
then
$$
\lim_{n \to +\infty}\int_{\mathbb{R}^{2}}F(u_n)\,dx=\lim_{n \to +\infty}\int_{\mathbb{R}^{2}}f(u_n)u_n\,dx=0.
$$
\end{lem}

As a consequence of the above lemma, we have the following result
\begin{cor} \label{sequencia}  Let $(u_n) \subset H^{1}(\mathbb{R}^{2})$ be a  $(PS)_c$ sequence for $I_\epsilon$ with $c>0$, $u_n \rightharpoonup 0$ and 
$$
\limsup_{n \to +\infty} \|u_n \|^{2}  < 1.
$$
Then, there exists $(z_n) \subset \mathbb{R}^{2}$ with $|z_n| \to +\infty$ such that
$$
v_n=u_n(\cdot+z_n) \rightharpoonup v \not=0 \quad \mbox{in} \quad H^{1}(\mathbb{R}^{2}).
$$
\end{cor}
\noindent {\bf Proof.}  We begin claiming that for any $R>0$, 
$$
\liminf_{n \to +\infty}\sup_{z \in \mathbb{R}^{2}}\int_{B_R(z)}|u_n|^{2}\,dx>0.
$$
Otherwise, there is $R>0$ and subsequence of $(u_n)$, still denoted by $(u_n)$ such that   
$$
\liminf_{n \to +\infty}\sup_{z \in \mathbb{R}^{2}}\int_{B_R(z)}|u_n|^{2}\,dx=0.
$$ 
Hence, by Lemma \ref{lions}, 
$$
\lim_{n \to +\infty}\int_{\mathbb{R}^{2}}f(u_n)u_n\,dx=0.
$$
The last limit combined with $I'_\epsilon(u_n)u_n=o_n(1)$ gives 
$$
u_n \to 0 \quad \mbox{in} \quad H^{1}(\mathbb{R}^{2}). 
$$ 
Then,
$$
I_\epsilon(u_n) \to 0,
$$
which is a contradiction, because by hypotheses $I_\epsilon(u_n) \to c>0$. Thereby, for each $R>0$,  there are $(z_n) \subset \mathbb{R}^{2}$, $\tau>0$ and a subsequence of $(u_n)$, still denoted by itself, such that
\begin{equation} \label{L 1}
\int_{B_R(z_n)}|u_n|^{2}\,dx \geq \tau.
\end{equation} 
Setting $v_n=u_n(\cdot+z_n)$, we have that $(v_n)$ is bounded in $H^{1}(\mathbb{R}^{2})$. Thus, for some subsequence of $(v_n)$, still denoted by $(v_n)$, 
there is $v \in H^{1}(\mathbb{R}^{2})$ such that
\begin{equation} \label{L2}
v_n \rightharpoonup v \quad \mbox{in} \quad H^{1}(\mathbb{R}^{2}).
\end{equation} \label{L3}
From (\ref{L 1}) and (\ref{L2}),
$$
\int_{B_R(0)}|v|^{2}\,dx \geq \tau,
$$
showing that $v \not= 0$. Moreover, (\ref{L 1}) guarantees that $|z_n| \to +\infty$, because $u_n \rightharpoonup 0$ in $H^{1}(\mathbb{R}^{2})$.
\fim

\vspace{0.5 cm}

The next lemma will be used to show the $(PS)$ condition at some level for $I_\epsilon$.

\begin{lem} \label{decomposicao}
	Let $(u_n) \subset H^{1}(\mathbb{R}^{2})$ be a sequence  verifying
\begin{equation} \label{H}
\limsup_{n \to +\infty} \|u_n \|^{2}  < \frac{1}{4},
\end{equation}
and $v_n=u_n-u$. Then,\\
\noindent a) \, $\displaystyle \lim_{n \to +\infty}\int_{\mathbb{R}^{2}}|F(v_n)-F(u_n)-F(u)|\,dx=0.$\\
\noindent b) \, There is $r>1$ with $ r \approx 1 $ such that
$$
\lim_{n \to +\infty}\int_{\mathbb{R}^{2}}|f(v_n)-f(u_n)+f(u)|^{r}\,dx=0.
$$ 
\end{lem}
\noindent {\bf Proof.}
From definition of $v_n$,
\begin{equation} \label{v1}
\limsup_{n \to \infty}\|v_n\|^{2}<\frac{1}{4}.
\end{equation}
Then, if $w_n=|v_n|+|u|$,  
$$
\limsup_{n \to +\infty}\|w_n\|^{2}<1,
$$
implying that there are $m \in (0,1)$ and $n_0 \in \mathbb{N}$ satisfying
\begin{equation} \label{W1}
\|w_n\|^{2}\leq m, \quad \forall n \geq n_0.
\end{equation}
By using Mean Value Theorem and $(f_1)$, 
$$
|F(v_n)-F(u_n)|\leq Ce^{4 \pi|w_n|^{2}}|u|, \quad \forall n \in \mathbb{N}.
$$
Thus, for any  $\beta >1$ and  $R>0$, 
$$
\int_{|x|>R}|F(v_n)-F(u_n)|\,dx \leq  C \int_{|x|>R}h_n|u|\,dx,
$$
where $h_n=e^{4 \beta \pi  |w_n|^{2}}-1$. Fixing $\beta \approx 1$ such that $\beta m <1$, by (\ref{W1}) and Lemma \ref{alphat11}, there is $q>1$ with $q \approx 1$ with
$$
h_n \in L^{q}(\mathbb{R}^{2}) \quad \mbox{and} \quad \sup_{n \in \mathbb{N}}|h_n|_q<+\infty.
$$
Now, applying H\"{o}lder's inequality, 
$$
\int_{|x|>R}h_n|u|\,dx \leq C |u|_{L^{q'}(|x|>R)} \quad \forall n \in \mathbb{N},
$$
where $\frac{1}{q}+\frac{1}{q'}=1$. From this, given $\delta>0$, we can choose $R>0$ large enough satisfying $$
\int_{|x|>R}h_n|u|\,dx < \frac{\delta}{4}, \quad \forall n \in \mathbb{N}.
$$
Hence, increasing $R$ if necessary, 
\begin{equation} \label{EST1}
\int_{|x|>R}|F(v_n)-F(u_n)-F(u)|\,dx < \frac{\delta}{2} \quad \forall n \in \mathbb{N}.
\end{equation}
On the other hand, from Corollary \ref{Convergencia em limitados} ,
\begin{equation} \label{EST2}
\lim_{n \to +\infty}\int_{|x|\leq R}|F(v_n)-F(u_n)-F(u)|\,dx=0.
\end{equation}
From (\ref{EST1}) and (\ref{EST2})
$$
\limsup_{n \to +\infty}\int_{\mathbb{R}^{2}}|F(v_n)-F(u_n)-F(u)|\,dx \leq \frac{\delta}{2}. 
$$
As $\delta$ is arbitrary, we can conclude that
$$
\lim_{n \to +\infty}\int_{\mathbb{R}^{2}}|F(v_n)-F(u_n)-F(u)|\,dx=0,
$$
showing a). To proof b), we use a similar approach. Applying the Mean Value Theorem together with $(f_1)$, we find  
$$
|f(v_n)-f(u_n)|\leq Ce^{4 \pi|w_n|^{2}}|u| \quad \forall n \in \mathbb{N}.
$$
Then, for all $R,r>0$  
$$
\int_{|x|>R}|f(v_n)-f(u_n)|^{r}\,dx \leq  C \int_{|x|>R}g_n|u|^{r}\,dx,
$$
where $g_n=e^{4 r \pi|w_n|^{2}}-1$. By (\ref{W1}), there are $s,r>1$ with $q,r \approx 1$ such that 
$$
g_n \in L^{s}(\mathbb{R}^{2}) \quad \mbox{and} \quad \sup_{n \in \mathbb{N}}|g_n|_s<+\infty.
$$
From H\"{o}lder's inequality, 
$$
\int_{|x|>R}g_n|u|^{r}\,dx \leq C |u|^{r}_{L^{s'}(|x|>R)} \quad \forall n \in \mathbb{N},
$$
where $\frac{1}{s}+\frac{1}{s'}=1$. Then, given $\delta>0$, we can fix $R>0$ large enough such that
$$
\int_{|x|>R}h_n|u|\,dx < \frac{\delta}{4} \quad \forall n \in \mathbb{N}.
$$
Again, increasing $R$ if necessary, 
\begin{equation} \label{EST3}
\int_{|x|>R}|f(v_n)-f(u_n)-f(u)|^{r}\,dx < \frac{\delta}{2} \quad \forall n \in \mathbb{N}.
\end{equation}
Now, repeating the same arguments explored in the proof of Corollary \ref{Convergencia em limitados}, we derive that
\begin{equation} \label{EST4}
\lim_{n \to +\infty}\int_{|x|\leq R}|f(v_n)-f(u_n)-f(u)|^{r}\,dx=0.
\end{equation}
From (\ref{EST3}) and (\ref{EST4}),
$$
\limsup_{n \to +\infty}\int_{\mathbb{R}^{2}}|f(v_n)-f(u_n)-f(u)|^{r}\,dx\leq \frac{\delta}{2}.
$$
As $\delta >0$ is arbitrary, we deduce that
$$
\lim_{n \to +\infty}\int_{\mathbb{R}^{2}}|f(v_n)-f(u_n)-f(u)|^{r}\,dx=0,
$$
showing b). \fim

\vspace{0.5 cm}

The lemma below brings an important estimate from above involving the mountain pass level $m(c_0)$, which will permit to use the Trundiger-Moser inequality found in \cite{Cao}

\begin{lem} \label{ESTIMATIVA SUPERIOR} Let $(u_n) \subset H^{1}(\mathbb{R}^{2})$ be a sequence $(PS)_c$ for $I_\epsilon$ with $c \in (0,2m(c_0))$. Then, 
	$$
	\limsup_{n \to +\infty}\|u_n\|^{2}< \frac{1}{4}.
	$$	
\end{lem}
\noindent {\bf Proof.}
Consider $w \in H^{1}(\mathbb{R}^{2})$ verifying 
$$
J_{\infty}(w)=\beta_p \quad \mbox{and} \quad J_{\infty}'(w)=0. 
$$
By characterization of $m(c_0)$, 
$$
m(c_0) \leq \max_{t \geq 0}I_\epsilon(tw).
$$
Consequently, by $(f_5)$, 
$$
m(c_0) \leq \max_{t \geq 0}\left\{ \frac{t^{2}}{2}\int_{\mathbb{R}^{2}}(|\nabla w|^{2}+|V|_{\infty}|w|^{2})\,dx-\frac{C_pt^{p}}{p}\int_{\mathbb{R}^{2}}|w|^{p}\,dx\right\},
$$
implying that 
$$
m(c_0) \leq C_p^{\frac{2}{2-p}}\beta_p.
$$
Gathering $I_\epsilon(u_n)-\frac{1}{\theta}I_\epsilon'(u_n)u_n=c+o_n(1)$ with $(f_3)$, we find 
$$ 
\left( \frac{1}{2}-\frac{1}{\theta}\right)\int_{\mathbb{R}^{2}}(|\nabla u_n|^{2}+V(x)|u_n|^{2})\,dx \leq c+o_n(1),
$$
from where it follows that
$$
\min\{1,c_0\} \left( \frac{1}{2}-\frac{1}{\theta}\right)\limsup_{n \to +\infty}\|u_n\|^{2} \leq c \leq 2m(c_0) \leq 2 C_p^{\frac{2}{2-p}}\beta_p .
$$
Since 
$$
C_p>\left[8\beta_p\left(\frac{2\theta}{\theta-2}\right)\frac{1}{\min\{1,c_0\}}\right]^{(p-2)/2},
$$
we obtain
$$
\limsup_{n \to +\infty}\|u_n\|^{2} < \frac{1}{4}.
$$
\fim

After the previous technical lemmas, we are ready to study the $(PS)$ condition for $I_\epsilon$.

\begin{lem} \label{L1} Under the hypotheses $(V_1)-(V_4)$, for each $\sigma >0$, there is $\epsilon_0=\epsilon_0(\sigma)>0$, such that  $I_\epsilon$ satisfies the $(PS)_c$ condition for all $c \in (m(c_0)+\sigma,2m(c_0)-\sigma)$, for all $\epsilon \in (0, \epsilon_0)$. 
\end{lem}

\noindent {\bf Proof.}  We will prove the lemma arguing by contradiction, by supposing that there are $\sigma >0$ and $\epsilon_n \to 0$, such that $I_{\epsilon_n}$ does not satisfy the $(PS)$ condition.  

Thereby, there is $c_n \in  (m(c_0)+\sigma,2m(c_0)-\sigma)$ such that  $I_{\epsilon_n}$ does not satisfy the $(PS)_{c_n}$ condition. Then, there  is a sequence $(u^{n}_m)$ such that
\begin{equation} \label{E1}
\lim_{m \to +\infty}I_{\epsilon_n}(u^{n}_m)=c_n \quad \mbox{and} \quad \lim_{m \to +\infty}I'_{\epsilon_n}(u^{n}_m)=0,
\end{equation}
with
\begin{equation} \label{E2}
 u^{n}_m \rightharpoonup u_n \quad \mbox{in} \quad H^{1}(\mathbb{R}^{2}) \quad \mbox{but} \quad u^{n}_m \not\to u_n \quad \mbox{in} \quad H^{1}(\mathbb{R}^{2}).
\end{equation}
By Lemma \ref{ESTIMATIVA SUPERIOR}, 
$$
\limsup_{m \to +\infty}\|u^{n}_m\|^{2} < \frac{1}{4}. 
$$
Then, for $v^{n}_m=u^{n}_m-u_n$, the Lemma \ref{decomposicao} yields 
$$
I_{\epsilon_n}(u^{n}_m)=I_{\epsilon_n}(u_n)+I_{\epsilon_n}(v^{n}_m)+o_m(1) \quad \mbox{and} \quad I'_{\epsilon_n}(v^{n}_m)=o_m(1).
$$

\begin{claim} \label{CorLions} There is $\delta>0$, such that
$$
\liminf_{m \to +\infty}\sup_{y \in \mathbb{R}^{2}}\int_{B_R(y)}|v^{n}_{m}|^{2}\,dx \geq \delta, \quad \forall n \in \mathbb{N}.
$$	
\end{claim}
Indeed, if the claim does not hold, there is $(n_j) \subset \mathbb{N}$ satisfying
$$
\liminf_{m \to +\infty}\sup_{y \in \mathbb{R}^{2}}\int_{B_R(y)}|v^{n_j}_{m}|^{2}\,dx \leq \frac{1}{j}, \quad \forall j \in \mathbb{N}.
$$
Then, repeating the arguments found in \cite[Lemma 1.21]{Willem}, se see that
\begin{equation} \label{Lq}
\limsup_{m \to +\infty}|v^{n_j}_{m}|_q=o_j(1), \quad \forall q \in (2,+\infty).
\end{equation}
As
$$
\limsup_{m \to +\infty}\|v^{n_j}_m\|^{2} < \frac{1}{4}, 
$$
it follows that
$$
\limsup_{m \to +\infty}\int_{\mathbb{R}^{2}}f(v^{n_j}_{m})(v^{n_j}_{m})\,dx=o_j(1).
$$
The above estimate, (\ref{Lq}) and $I'_{\epsilon_{n_j}}(v^{n_j}_{m})(v^{n_j}_{m})=o_{m}(1)$ combine to give
$$
\limsup_{m \to +\infty}\|v^{n_j}_m\|^{2}=o_j(1).
$$
Since $u^{n_j}_m \not\to u_{n_j}$ in $H^{1}(\mathbb{R}^{2})$, we know that 
$$
\liminf_{m \to +\infty}\|v^{n_j}_m\|^{2}>0.
$$
Therefore, without loss of generality, we can assume that $(v^{n_j}_m) \subset H^{1}(\mathbb{R}^{2}) \setminus \{0\}$. Thereby, there is $t^{n_j}_m \in (0,+\infty)$ such that
$$
t^{n_j}_m v^{n_j}_m \in \mathcal{N}_{\epsilon_{n_j}}.
$$
Using $(f_4)$, it is possible to prove that
$$
\lim_{m \to +\infty}t^{n_{j}}_m=1
$$
and
$$
\lim_{m \to +\infty}I_{\epsilon_{n_j}}(t^{n_{j}}_m v^{n_j}_m)= \lim_{m \to +\infty}I_{\epsilon_{n_j}}(v^{n_j}_m).
$$ 
From the above informations, there is $r^{n_j}_m \in (0,1)$ such that
$$
r^{n_j}_m (t^{n_j}_m v^{n_j}_m) \in \mathcal{M}_{c_0}.
$$
Hence,
$$
m(c_0) \leq \limsup_{m \to +\infty}J_{c_0}(r^{n_j}_m (t^{n_j}_m v^{n_j}_m))\leq \limsup_{m \to +\infty}I_{\epsilon_{n_j}}(t^{n_j}_m v^{n_j}_m) =\limsup_{m \to +\infty}I_{\epsilon_{n_j}}(v^{n_j}_m)\leq \frac{(1+|V|_{\infty})}{2} \limsup_{m \to +\infty}\|v^{n_j}_m\|^{2},
$$
that is,
$$
0<m(c_0) \leq o_j(1),
$$
which is a contradiction.

From the above study, for each $m \in \mathbb{N}$, there is $m_n \in \mathbb{N}$ such that
$$
\int_{B_R(z^{n}_{}m_n)}|u^{n}_{m_n}|^{2}\,dx \geq \frac{\delta}{2}, \quad |\epsilon_n z^{n}_{m_n}| \geq n, \quad \|I'_{\epsilon_n}(u^{n}_{m_n})\|\leq \frac{1}{n} \quad \mbox{and} \quad |I_{\epsilon_n}(u^{n}_{m_n})-c_n|\leq \frac{1}{n}.
$$

In what follows, we denote by $(z_n)$ and $(u_n)$ the sequences $(z^{n}_{m_n})$ and $(u^ {n}_{m_n})$ respectively.  Then,
$$
\int_{B_R(z_n)}|u_{n}|^{2}\,dx \geq \frac{\delta}{2}, \quad |\epsilon_n z_{n}| \geq n, \quad \|I'_{\epsilon_n}(u_n)\|\leq \frac{1}{n} \quad \mbox{and} \quad |I_{\epsilon_n}(u_n)-c_n|\leq \frac{1}{n}.
$$

\begin{claim} \label{limite fraco}  $(u_n)$ is a bounded sequence in $H^{1}(\mathbb{R}^{2}).$ Moreover, for some subsequence, $u_n \rightharpoonup 0$ in $H^{1}(\mathbb{R}^{2})$. 
	
\end{claim}

In fact, the boundedness of $(u_n)$ follows by standard arguments. Then, for some subsequence, there is $u \in H^{1}(\mathbb{R}^{2})$ such that
$$
u_n \rightharpoonup u \quad \mbox{in} \quad H^{1}(\mathbb{R}^{2}). 
$$
Supposing by contradiction that $u \not= 0$, the limit $ \|I'_{\epsilon_n}(u_n)\|\to 0$ together with $(V_2)$ yield $u$ is a nontrivial solution of the problem
$$
\Delta{u}-V(0)u+f(u)=0 \quad \mbox{in} \quad \mathbb{R}^{2}.
$$
Then,  combining the definition of $m(V(0))$  with  $(V_4)$, we get   
$$
J_{V(0)}(u) \geq m(V(0)) \geq 2m(c_0).
$$
On the other hand,  the Fatous' lemma loads to
$$
J_{V(0)}(u) \leq \liminf_{n \to +\infty}[I_{\epsilon_n}(u_n)-\frac{1}{\theta}I'_{\epsilon_n}(u_n)]=\liminf_{n \to +\infty}I_{\epsilon_n}(u_n)=\liminf_{n \to +\infty}c_n \leq 2m(c_0)-\sigma,
$$
obtaining a contradiction. From this, the proof of Claim \ref{limite fraco} is finished.

Considering $w_n=u_n(\cdot+z_n)$, we have that $(w_n)$ is bounded in $H^{1}(\mathbb{R}^{2})$. Therefore, there is $w \in H^{1}(\mathbb{R}^{2})$ such that
$$
w_n \rightharpoonup w \quad \mbox{in} \quad H^{1}(\mathbb{R}^{2})
$$
and
$$
\int_{B_R(0)}|w|^{2}\,dx=\liminf_{n \to +\infty}\int_{B_R(0)}|w_n|^{2}\,dx=\liminf_{n \to +\infty}\int_{B_R(z_n)}|u_n|^{2}\,dx \geq \frac{\delta}{2},
$$
showing that $w \not= 0$.

Now, for each $\phi \in H^{1}(\mathbb{R}^{2})$, we have the equality below
$$
\int_{\mathbb{R}^{2}}\nabla w_n \nabla \phi \, dx + \int_{\mathbb{R}^{2}}V(\epsilon_n z_n+\epsilon z)w_n \phi \, dx - \int_{\mathbb{R}^{2}}f(w_n)\phi \, dx = o_n(1)\|\phi\|,
$$
which implies that $w$ is a nontrivial solution of the problem
\begin{equation} \label{equacao}
\Delta{u}-\alpha_1u+f(u)=0 \quad \mbox{in} \quad \mathbb{R}^{2},
\end{equation}
where $\alpha_1=\displaystyle \lim_{n \to +\infty}V(\epsilon_n z_n)$. Thus, by regularity theory, $w \in C^{2}(\mathbb{R}^{2}) \cap H^{2}(\mathbb{R}^{2})$. 

For each $k \in \mathbb{N}$, there is $\phi_k \in C^{\infty}_{0}(\mathbb{R}^{2})$ such that
$$
\|\phi_k -w\| \to 0 \quad \mbox{as} \quad k \to +\infty,
$$
that is,
$$
\|\phi_k -w\|=o_k(1).
$$

Using $\frac{\partial \phi_k}{\partial x_i}$ as a test function, we get
$$
\int_{\mathbb{R}^{2}}\nabla w_n \nabla \frac{\partial \phi_k}{\partial x_i}\,dx+ \int_{\mathbb{R}^{2}}V(z+z_n)w_n\frac{\partial \phi_k}{\partial x_i}\,dx-
\int_{\mathbb{R}^{2}}f(w_n)\frac{\partial \phi_k}{\partial x_i}\,dx=o_n(1).
$$
Now, using well known arguments,
$$
\int_{\mathbb{R}^{2}}\nabla w_n \nabla \frac{\partial \phi_k}{\partial x_i}\,dx=\int_{\mathbb{R}^{2}}\nabla w \nabla \frac{\partial \phi_k}{\partial x_i}\,dx+o_n(1)
$$
and
$$
\int_{\mathbb{R}^{2}}f(w_n)\frac{\partial \phi_k}{\partial x_i}\,dx=\int_{\mathbb{R}^{2}}f(w)\frac{\partial \phi_k}{\partial x_i}\,dx+o_n(1).
$$
Gathering the above limit with (\ref{equacao}), we deduce that  
$$
\limsup_{n \to +\infty}\left|\int_{\mathbb{R}^{2}}(V(\epsilon_n z_n+\epsilon_n z)-V(\epsilon_n z_n))w_n\frac{\partial \phi_k }{\partial x_i}\right|\,dx=0.
$$
As $\phi_k$ has compact support, the above limit gives 
$$
\limsup_{n \to +\infty}\left|\int_{\mathbb{R}^{2}}(V(\epsilon_n z_n+\epsilon_n z)-V(z_n))w\frac{\partial \phi_k}{\partial x_i}\right|\,dx=0.
$$
Now, recalling that $\frac{\partial w}{\partial x_i} \in L^{2}(\mathbb{R}^{2})$, we have that $(\frac{\partial \phi_k}{\partial x_i})$ is bounded in $L^{2}(\mathbb{R}^{2}).$ Hence, 
$$
\limsup_{n \to +\infty}\left|\int_{\mathbb{R}^{2}}(V(\epsilon_n z_n+\epsilon_n z)-V(z_n))\phi_k\frac{\partial \phi_k}{\partial x_i}\right|\,dx=o_k(1),
$$
and so,
$$
\limsup_{n \to +\infty}\left|\frac{1}{2}\int_{\mathbb{R}^{2}}(V(\epsilon_n z_n+\epsilon_n z)-V(z_n))\frac{\partial (\phi_k^{2})}{\partial x_i}\right|=o_k(1).
$$
Using Green's Theorem together with the fact that $\phi_k$ has compact support, we find the limit below
$$
\limsup_{n \to +\infty}\left|\int_{\mathbb{R}^{2}}\frac{\partial V}{\partial x_i}(\epsilon_n z_n+\epsilon_n z) \, \phi_k^{2}\right|\,dx=o_k(1),
$$
which combined with $(V_2)$ gives 
$$
\limsup_{n \to +\infty}\left|\frac{\partial V}{\partial x_i}(\epsilon_n z_n)\int_{\mathbb{R}^{2}}|\phi_k|^{2}\right|\,dx=o_k(1).
$$
As 
$$
\int_{\mathbb{R}^{2}}|\phi_k|^{2}\,dx \to \int_{\mathbb{R}^{2}}|w|^{2}\,dx \quad \mbox{as} \quad k \to +\infty,
$$
it follows that
$$
\limsup_{n \to +\infty}\left|\frac{\partial V}{\partial x_i}(\epsilon_n z_n)\right|=o_k(1), \quad \forall i \in \{1,....,N\}.
$$
Since $k$ is arbitrary, we derive that 
$$
\nabla V(\epsilon_n z_n) \to 0 \quad \mbox{as} \quad n \to \infty.
$$
Therefore,  $(\epsilon_n z_n)$ is a $(PS)_{\alpha_1}$ sequence for $V$, which is an absurd, because by hypotheses $V$ satisfies the $(PS)$ condition and $(z_n)$ does not have any convergent subsequence in $\mathbb{R}^{2}$. 
\fim

\vspace{0.5 cm}

Hereafter, we denote by ${\mathcal N}_{\epsilon}$ the Nehari Manifold associated with $I_\epsilon$, that is,
$$ 
 {\mathcal N}_{\epsilon}=\left\{ u \in H^{1}(\mathbb{R}^{2}) \setminus \{0\}\,:\, I'_{\epsilon}(u)u=0    \right\}.
$$

\begin{lem}\label{PS em Nehari} For each $\sigma>0$, there is $\epsilon_0=\epsilon_0(\sigma)>0$ such that for $\epsilon \in (0, \epsilon_0)$, the functional $I_\epsilon$ restrict to ${\mathcal N}_{\epsilon}$ satisfies the $(PS)_c$ condition for all $c \in (m(c_0)+\sigma,2m(c_0)-\sigma)$. 
\end{lem}
\noindent {\bf Proof.} Let $(u_n)$ be a $(PS)$ sequence for
$I_{\epsilon}$ constrained to $\mathcal{N}_{\epsilon}$. Then
$I_{\epsilon}(u_{n})\rightarrow c$ and
\begin{eqnarray}\label{contra}
I'_{\epsilon}(u_{n}) = \theta_{n} G_{\epsilon}'(u_{n}) + o_{n}(1),
\end{eqnarray}
for some $(\theta_{n}) \subset \mathbb{R}$, where
$G_{\epsilon}:H^{1}(\mathbb{R}^{2}) \rightarrow \mathbb{R}$ is
given by
\begin{eqnarray*}
	G_{\epsilon}(v) := \displaystyle\int_{\mathbb{R}^{2}}(|\nabla v|^{2}+V(\epsilon x)|v|^{2}) \,dx
	-\displaystyle\int_{\mathbb{R}^{2}}f(v)v \,dx.
\end{eqnarray*}
We recall that $G'_{\epsilon}(u_{n})u_{n}\leq 0$. Moreover, standard
arguments show that $(u_n)$ is bounded. Thus, up to a subsequence,
$G_{\epsilon}'(u_{n}) u_{n} \rightarrow l\leq 0$. If $l \neq 0$, we infer
from (\ref{contra}) that $\theta_{n}=o_{n}(1)$. In this case, we can
use (\ref{contra}) again to conclude that $(u_{n})$ is a $(PS)_{c}$
sequence for $I_{\epsilon}$ in $H^{1}(\mathbb{R}^{2})$, and so, $(u_{n})$ has a strongly convergent subsequence. If
$l=0$, it follows that 
$$ 
\displaystyle\int_{\mathbb{R}^{2}}(f'(u_n)u_n^{2}-f(u_n)u_n) \ dx\rightarrow 0.
$$
Using $(f_4)$, we know that
\begin{equation} \label{Desig}
f'(t)t^{2}-f(t)t>0 \quad \forall t \in \mathbb{R} \setminus \{0\}.
\end{equation}
Thereby, if $u \in H^{1}(\mathbb{R}^{2})$ is the weak limit of $(u_n)$, the Fatous' Lemma combined with the last limit gives
$$
\int_{\mathbb{R}^{2}}(f'(u)u^{2}-f(u)u) \ dx=0.
$$
Then, by (\ref{Desig}), $u=0$. Applying Corollary \ref{sequencia}, there is $(y_n) \subset \mathbb{R}^{2}$ with $|y_n| \to +\infty$ such that
$$
v_n=u_n(\cdot +y_n) \rightharpoonup v \not= 0 \quad \mbox{in} \quad H^{1}(\mathbb{R}^{2}).
$$
By change variable, we have that
$$
\displaystyle\int_{\mathbb{R}^{2}}(f'(v_n)v_n^{2}-f(v_n)v_n)\,dx=\displaystyle\int_{\mathbb{R}^{2}}(f'(u_n)u_n^{2}-f(u_n)u_n) \ dx\rightarrow 0.
$$
Applying again Fatous's Lemma, we get
$$
\int_{\mathbb{R}^{2}}(f'(v)v^{2}-f(v)v) \ dx=0,
$$
which is an absurd, because being  $v \not=0$, the inequality (\ref{Desig}) loads to
$$
\int_{\mathbb{R}^{2}}(f'(v)v^{2}-f(v)v) \ dx>0,
$$
finishing the proof of the lemma.

\fim

\begin{cor} \label{ponto critico}
If $u \in H^{1}(\mathbb{R}^{2})$ is a critical point of $I_\epsilon$ restrict to ${\mathcal N}_{\epsilon}$, then $u$ is a critical point of  $I_\epsilon$ in $H^{1}(\mathbb{R}^{2})$. 
\end{cor}
\noindent {\bf Proof.} The corollary follows adapting the arguments explored in the proof of Lemma \ref{PS em Nehari}. \fim

The next lemma will be crucial in our study to show an estimate from below involving  a special minimax level,  which will be considered later on. 
\begin{lem} \label{compacidade} Let $\epsilon_n \to 0$ and $(u_n) \subset \mathcal{N}_{\epsilon_n}$ such that $I_{\epsilon_n}(u_n) \to m(c_0)$. Then, there is $(z_n) \subset \mathbb{R}^{2}$ with $|z_n| \to +\infty$ and $u_1 \in H^{1}(\mathbb{R}^{2}) \setminus \{0\}$ such that
$$
u_n(\cdot +z_n) \to u_1 \quad \mbox{in} \quad H^{1}(\mathbb{R}^{2}).
$$
Moreover, $\displaystyle\liminf_{n \to +\infty}|\epsilon_n z_n|>0$.	
\end{lem}
\noindent {\bf Proof.} Since $u_n \in \mathcal{N}_{\epsilon_n}$, we have that
$J_{c_0}'(u_n)u_n<0$ and $J_{c_0}(u) \leq I_{\epsilon_n}(u)$ for all $u \in H^{1}(\mathbb{R}^{2})$ and $n \in \mathbb{N}$. From this, there is $t_n \in (0, 1)$ such that   
$$
(t_n u_n) \subset  \mathcal{M}_{c_0} \quad \mbox{and} \quad J_{c_0}(t_n u_n)\to m(c_0).
$$
Now, using \cite[Lemma 12]{AlvesGio}, there are $(z_n) \subset \mathbb{R}^{2}$, $u_1 \in H^{1}(\mathbb{R}^{2}) \setminus \{0\}$, and a subsequence of $(u_n)$, still denoted by $(u_n)$,  verifying
$$
u_n(\cdot+z_n) \to u_1 \quad \mbox{in} \quad H^{1}(\mathbb{R}^{2}).
$$
\begin{claim} \label{zn} 
$\displaystyle \liminf_{n \to +\infty}|\epsilon_n z_n|>0 $.	
\end{claim}
Indeed, as $u_n \in \mathcal{N}_{\epsilon_n}$ for all $n \in \mathbb{N}$, the function $u_n^{1}=u_n(\cdot+z_n)$ must verify
\begin{equation} \label{Eq1}
\int_{\mathbb{R}^{2}}(|\nabla u_{n}^{1}|^{2}+V(\epsilon_n x+\epsilon_n z_n)|u_n^{1}|^{2})\,dx=\int_{\mathbb{R}^{2}}f(u_n^{1})u_n^{1}\,dx.
\end{equation}
Supposing by contradiction that for some subsequence
$$
\lim_{n \to +\infty}\epsilon_n z_n=0,
$$
taking the limit of $n \to +\infty$ in (\ref{Eq1}), we derive the equality below
$$
\int_{\mathbb{R}^{2}}(|\nabla u_{1}|^{2}+V(0)|u_{1}|^{2})\,dx=\int_{\mathbb{R}^{2}}f(u_{1})u_{1}\,dx,
$$
showing that $u_1 \in \mathcal{M}_{V(0)}$. Thereby, 
\begin{equation} \label{EQm0}
J_{V(0)}(u_{1}) \geq m(V(0)) > m(c_0).
\end{equation}
On the other hand, by $(V_2)$,
$$
I_{\epsilon_n}(u_n) \to J_{V(0)}(u_{1}),
$$ 
loading to
\begin{equation} \label{EQm}
m(c_0)=J_{V(0)}(u_{1}).
\end{equation}
From (\ref{EQm0}) and (\ref{EQm}), we find a contradiction, finishing the proof.  
\fim	
	
\section{A special minimax level}

In order to prove the Theorem \ref{T1}, we will consider  a special minimax level. To do that, we begin fixing the barycenter function  by
$$
\beta(u)= \frac{\displaystyle  \int_{\mathbb{R}^{2}}\frac{x}{|x|}|u|^{2}\,dx}{\displaystyle \int_{\mathbb{R}^{2}}|u|^{2}\,dx}, \quad \forall u \in H^{1}(\mathbb{R}^{2}) \setminus \{0\}.
$$
In what follows, $u_0$ denotes a radial positive ground state solution for $J_{c_0}$, that is, 
$$
J_{c_0}(u_0)=m(c_0) \quad \mbox{and} \quad J'_{c_0}(u_0)=0.
$$ 
For each $z \in \mathbb{R}^{2}$ and $\epsilon >0$, we set the function
$$
\phi_{\epsilon, z}(x)=t_{\epsilon,z}u_0\left( x - \frac{z}{\epsilon} \right),
$$
where $t_{\epsilon,z}>0$ is such that $\phi_{\epsilon, z} \in {\mathcal N}_{\epsilon}$. From definition of $\beta$, we have the following result

\begin{lem} \label{B1} 
	For each  $r>0$, $\displaystyle \lim_{\epsilon \to 0}\left(\sup\left\{\left|\beta(\phi_{\epsilon,z})-\frac{z}{|z|} \right|\,:\,|z| \geq r \right\}\right)=0$.
\end{lem}
\noindent {\bf Proof.} The proof follows showing that for any $(z_n) \subset \mathbb{R}^{2}$ with $|z_n| \geq r$ and $\epsilon_n \to 0$, we have that 
$$
\left|\beta(\phi_{\epsilon_n,z_n})-\frac{z_n}{|z_n|} \right| \to 0 \quad \mbox{as} \quad n \to +\infty.
$$
By change variable, 
$$
\left|\beta(\phi_{\epsilon_n,z_n})-\frac{z_n}{|z_n|} \right|=\frac{\displaystyle \int_{\mathbb{R}^{2}}\left|\frac{\epsilon_n x +z_n}{|\epsilon_n x +z_n|}-\frac{z_n}{|z_n|}\right||u_0(x)|^{2}\,dx}{\displaystyle \int_{\mathbb{R}^{2}}|u_0|^{2}\, dx}.
$$
As for each $x \in \mathbb{R}^{2}$,
$$
\left|\frac{\epsilon_n x +z_n}{|\epsilon_n x +z_n|}-\frac{z_n}{|z_n|}\right| \to 0 \quad \mbox{as} \quad n \to +\infty, 
$$
the Lebesgue's Theorem ensures that  
$$
\int_{\mathbb{R}^{2}}\left| \frac{\epsilon_n x +z_n}{|\epsilon_n x +z_n|}-\frac{z_n}{|z_n|}\right||u_0(x)|^{2}\,dx \to 0,
$$
showing the lemma.  \fim

\vspace{0.5 cm}

As a by product of the arguments explored in the proof of the last lemma, we have the following corollary 
\begin{cor} \label{cor1} 
	Fixed $r>0$, there is $\epsilon_0>0$ such that
	$$
	(\beta(\phi_{\epsilon,z}),z)>0, \quad \forall |z| \geq r \quad \mbox{and} \quad \forall \epsilon \in (0,\epsilon_0).
	$$
\end{cor}

\vspace{0.5 cm}

In the sequel,  we define the set 
$$
{\mathcal B}_{\epsilon}=\{u \in {\mathcal N}_{\epsilon}\,:\, \beta(u) \in Y  \}.
$$
Note that ${\mathcal B}_{\epsilon} \not=\emptyset$, because $\phi_{\epsilon,0}=0 \in Y,$ for all $\epsilon>0$. 
Associated with the above set, we consider the real number $D{_\epsilon}$ given by
$$
D_{\epsilon}=\inf_{u \in \mathcal{B}_{\epsilon}}I_{\epsilon}(u).
$$

The next lemma establishes an important relation involving the levels $D_{\epsilon}$ and $m(c_0)$.
\begin{lem} \label{estimativas} \mbox{}\\
\noindent (a) \, There exist $\epsilon_0,\sigma>0$ such that
$$
D_{\epsilon} \geq m(c_0)+\sigma, \quad \forall \epsilon \in (0,\epsilon_0).
$$
\noindent (b) \, $ \displaystyle \limsup_{\epsilon \to 0}\left\{\sup_{x \in X}I_\epsilon(\phi_{\epsilon,x}) \right\} < 2m(c_0)-\sigma$. \\
\noindent (c) \, There exist $\epsilon_0,R>0$ such that
$$
I_\epsilon(\phi_{\epsilon,x}) \leq \frac{1}{2}(m(c_0)+D_\epsilon), \quad \forall \epsilon \in (0,\epsilon_0) \quad \mbox{and} \quad \forall x \in \partial B_R(0) \cap X.
$$ 
\end{lem}
\noindent {\bf Proof.} First we prove (a), arguing by contradiction. From definition of $D_{\epsilon} $, we know that
$$
D_{\epsilon}  \geq m(c_0), \quad \forall \epsilon >0. 
$$
Supposing by contradiction that the lemma does not hold, there exists $\epsilon_n \to 0$ verifying 
$$
D_{\epsilon_n} \to m(c_0).
$$
Hence, there is $u_n \in {\mathcal N}_{\epsilon_n}$, with $\beta(u_n)\in Y$, satisfying 
$$
I_{\epsilon_n}(u_n) \to m(c_0).
$$
Applying Lemma \ref{compacidade}, there are $u_1 \in H^{1}(\mathbb{R}^{2}) \setminus \{0\}$ and a sequence  $(z_n) \subset \mathbb{R}^{2}$  with $\displaystyle \liminf_{n \to +\infty}|\epsilon_n z_n|>0$ verifying 
$$
u_n(\cdot+z_n) \to u_1 \quad \mbox{in} \quad  H^{1}(\mathbb{R}^{2}),
$$
that is,
$$
u_n=u_1(\cdot-z_n)+w_n \quad \mbox{with} \quad w_n \to 0 \quad \mbox{in} \quad H^{1}(\mathbb{R}^{2}).
$$
The definition of $\beta$ and the same arguments explored in the proof of Lemma \ref{B1} combine to give 
$$
\beta(u_n)=\beta(u_1(\cdot-z_n))+o_n(1)=\frac{z_n}{|z_n|}+o_{n}(1).
$$
As $\beta(u_n) \in Y$, we infer that $\frac{z_n}{|z_n|} \in Y_\lambda$ for $n$ large enough. Consequently, $\epsilon_n z_n \in Y_\lambda$ for $n$ large enough, implying that
$$
\liminf_{n \to \infty}V(\epsilon_n z_n) > c_0.
$$
If $A=\displaystyle \liminf_{n \to \infty}V(\epsilon_n z_n)$, the last inequality together with Fatous's Lemma yields 
$$
m(c_0)=\liminf_{n \to \infty}J_{\epsilon_n}(u_n) \geq J_A(u_1) \geq  m(A)>m(c_0),
$$
which is an absurd.  Here, we have used the fact that $J'_A(u_1)u_1=0$ and $u_1 \not= 0$. \\

\noindent {\bf Proof of (b).} \, Using condition $(V_4)$ together with the fact that $u_0$ is a ground state solution associated with $I_{c_0}$, we deduce that 
$$
\limsup_{\epsilon \to 0}\left\{\sup_{x \in X}I_{\epsilon}(\phi_{\epsilon,x})\right\}\leq J_{c_0}(u_0)+\frac{3}{5}J_{c_0}(u_0)=m(c_0)+ \frac{3}{5}m(c_0)<2m(c_0) \quad \forall \epsilon \in (0, \epsilon_0).  
$$ 
\noindent {\bf Proof of (c).} \, From $(V_1)$, given $\delta >0$, there are $R,\epsilon_0>0$ such that
$$
\sup\{I_\epsilon(\phi_{\epsilon,x})\,:\, x \in \partial B_{R}(0) \cap X \}\leq m(c_0)+\delta, \quad \forall \epsilon \in (0, \epsilon_0). 
$$ 
Fixing $\delta = \frac{\sigma}{4}$, where $\sigma$ was given in (a), we  have that
$$
\sup\{I_\epsilon(\phi_{\epsilon,x})\,:\, x \in \partial B_{R}(0) \cap X \}\leq \frac{1}{2}\left(2m(c_0)+\frac{\sigma}{2}\right)< \frac{1}{2}(m(c_0)+D_\epsilon), \quad \forall \epsilon \in (0, \epsilon_0). 
$$

\fim

\vspace{0.5 cm}

Now, we are ready to show the minimax level. Define the map $\Phi_{\epsilon}:X \to H^{1}(\mathbb{R}^2)$ as  $\Phi_{\epsilon}=\phi_{\epsilon,x}$. In what follows, $P$ denotes the cone of nonnegative functions of $H^{1}(\mathbb{R}^{2})$ and we consider the set
$$
\Sigma=\{\Phi_\epsilon\,:\, x \in X, \, |x|\leq R \} \subset P,
$$
the class of functions
$$
\mathcal{H}=\left\{h \in C(P \cap \mathcal{N}_{\epsilon},P \cap \mathcal{N}_{\epsilon})\,:\,h(u)=u, \,\, \mbox{if} \,\, I_\epsilon(u) < \frac{1}{2}(m(c_0)+D_\epsilon) \right\} 
$$
and finally the class of sets
$$
\Gamma=\{A \subset P \cap \mathcal{N}_\epsilon\,:\, A=h(\Sigma), \, h \in \mathcal{H} \}.
$$
\begin{lem} \label{intersecao}
If $A \in \Gamma$, then $A \cap \mathcal{B}_\epsilon \not= \emptyset$ for all $\epsilon>0$.
\end{lem}
\noindent {\bf Proof.} It is enough to show that for all $h \in  \mathcal{H}$, there is $x_* \in X$ with $|x_*| \leq R$ such that
$$
\beta(h(\Phi_\epsilon(x_*))) \in Y.
$$
For each $h \in \mathcal{H}$, we set the function $g:\mathbb{R}^{2} \to \mathbb{R}^{2}$ given by
$$
g(x)=\beta(h(\Phi_\epsilon(x))) \quad \forall x \in \mathbb{R}^{2},
$$
and the homotopy $\mathcal{F}:[0,1] \times X \to X$ as
$$
\mathcal{F}(t,x)=tP_X(g(x))+(1-t)x,
$$
where $P_X$ is the projection onto $X$. By Corollary \ref{cor1} and Lemma \ref{B1}, fixed $R>0$ and $\epsilon>0$ small enough, we have that
$$
(\mathcal{F}(t,x),x)>0, \quad \forall (t,x)  \in [0,1] \times ( \partial B_R \cap X).
$$
Using the Topological degree, we derive
$$
d(g,B_R \cap X,0)=1,
$$
implying that there exists $x_* \in B_R \cap X$ such that  $ \beta(h(\Phi_\epsilon(x_*)))=0$.
\fim

Now, we define the min-max value
$$
C_\epsilon=\inf_{A \in \Gamma}\sup_{u \in A}I_{\epsilon}(u).
$$
From Lemmas \ref{estimativas} and \ref{intersecao},
\begin{equation} \label{NM1}
C_\epsilon \geq D_\epsilon \geq m(c_0)+\sigma,
\end{equation}
for $\epsilon$ is small enough. On the other hand, 
$$
C_\epsilon \leq \sup_{x \in X}I_\epsilon(\phi_{\epsilon,x}), \quad \forall \epsilon >0.
$$ 
Then, by Lemma \ref{estimativas}(b), if $\epsilon$ is small enough
\begin{equation} \label{NM2}
C_\epsilon \leq \sup_{x \in X}I_\epsilon(\phi_{\epsilon,x})< 2m(c_0)-\sigma.
\end{equation}
From (\ref{NM1}) and (\ref{NM2}), there is $\epsilon_0>0$ such that
\begin{equation} \label{NM3}
C_\epsilon \in (m(c_0)+\sigma,2m(c_0)-\sigma), \quad \forall \epsilon \in (0, \epsilon_0). 
\end{equation}

Now, we can use standard min-max arguments to conclude that $I_\epsilon$ has at least a critical point in $P \cap \mathcal{N}_{\epsilon}$ if $\epsilon$ is small enough.

\end{document}